\documentclass[10pt,letterpaper]{amsart}

\newcommand{\dealii}{\texttt{deal.II}}
\newcommand{\aspect}{\texttt{ASPECT}}

\usepackage{amsaddr}

\usepackage{amsmath}
\usepackage{amssymb}
\usepackage{epsfig}
\usepackage{graphics}
\usepackage{graphicx}
\usepackage{subfigure}
\usepackage{color}
\usepackage{multicol}
\usepackage{multimedia}
\usepackage{eso-pic} 
\usepackage{colortbl}
\usepackage{comment}
\usepackage{tabularx}

\usepackage{graphicx}              
\graphicspath{{../results/}}
\usepackage{caption}
\usepackage{subfigure}

\newcommand{\bvec}[1]{\bm{#1}}
\newcommand{\QQ}{\mathbb{Q}}

\usepackage{amsmath,bm}               
\usepackage{amsthm}                

\usepackage{mathtools}

\calclayout

\begin{document}

\title[Matrix-free Geometric Multigrid for a Stokes Problem]{Comparison Between Algebraic and Matrix-free Geometric Multigrid for a Stokes Problem on Adaptive Meshes with Variable Viscosity}

\author{Thomas C. Clevenger}
\address{{tcleven@clemson.edu}, Clemson University}
%
%
\author{Timo Heister}
\address{{heister@clemson.edu}, Clemson University}

%
%
%
%
%
%
%




\begin{abstract}
Problems arising in Earth's mantle convection involve finding the solution to Stokes systems with large viscosity contrasts. These systems contain localized features which, even with adaptive mesh refinement, result in linear systems that can be on the order of $10^9$ or more unknowns. One common approach for preconditioning to the velocity block of these systems is to apply an Algebraic Multigrid (AMG) v-cycle (as is done in the ASPECT software, for example), however, with AMG, robustness can be difficult with respect to problem size and number of parallel processes. Additionally, we see an increase in iteration counts with adaptive refinement when using AMG. In contrast, the Geometric Multigrid (GMG) method, by using information about the geometry of the problem, should offer a more robust option.

Here we present a matrix-free GMG v-cycle which works on adaptively refined, distributed meshes, and we will compare it against the current AMG preconditioner (Trilinos ML) used in the \aspect{}\cite{aspect1} software. We will demonstrate the robustness of GMG with respect to problem size and show scaling up to 114688 cores and $217$ billion unknowns. All computations are run using the open source, finite element library \dealii{}.\cite{dealII91}

\end{abstract}



\maketitle

\section{Introduction}
The major bottleneck of computations of processes in Earth's mantle convection is the solution of Stokes systems to solve for velocity and pressure. These systems often have large variation in their coefficients, as well as highly localized features requiring adaptive mesh refinement to yield high enough resolution while still being computationally feasible to solve. Even with adaptive refinement, there is a desire to solve on meshes containing well over a billion unknowns, which many current open-source codes are not equipped to handle.


When open-source codes like the mantle convection code \aspect{}\cite{aspect1} were first developed, the  state-of-the art solvers often relied on algebraic multigrid (AMG) methods when preconditioning on the velocity space of the Stokes finite element discretizations (see, e.g., the works of Geenen et al.\cite{geenen2009} and Kronbichler et al.\cite{kronbichler_heister_bangerth_2012}). While these algebraic methods can be very powerful for smaller problems, their performance tends to deteriorate with highly adaptive meshes and when distributing the problem over a large number of cores. In addition, these methods require the storing of matrices which greatly limits the size of the problem which can be solved, and can provide a major computational bottleneck as the access of matrix entries from main memory starts to become slower than the actual computation with those entries.\cite{kronbichler2016comparison} 

In response to these limitations, much of the research into solving the Stokes systems is now focused on using either geometric multigrid (GMG) or a hybrid multigrid (HMG) (where AMG is used on coarse meshes of the GMG method) for preconditioning on the velocity space. These methods have the advantage that, by using geometric information about the problem to construct the multigrid level hierarchy, there should be less deterioration with highly adaptive meshes. Also, using these geometric methods allow for the implementation of matrix-free operator evaluation which not only reduces the memory required for an application (allowing one to solve larger problems), but also can lead to faster matrix-vector products, particularly for higher order finite elements.\cite{kronbichler2012generic} These methods have lead to solutions of systems on the order of $10^{13}$ unknowns\cite{Gmeiner2015} (low order, globally refined tetrahedral mesh) and $10^{11}$ unknowns\cite{RudiMalossi2015} (higher order, adaptively refined hexahedral mesh).

Here we will present a comparison of a Stokes application using the matrix-free GMG method presented in the work of Clevenger et al.\cite{clevenger_preprint} and implemented by the authors into \aspect{} with the current AMG-based method described in the work of Kronbichler et al.\cite{kronbichler_heister_bangerth_2012} and Heister et al.\cite{heister_dannberg_gassmoeller_bangerth_2017} (with Trilinos ML AMG) widely used by the \aspect{} community. We will demonstrate the advantages of the geometric method over the algebraic method, as well as show weak and strong scalability of the geometric algorithm on the sinker benchmark used in the work of May et al.\cite{MAY2015496} and Rudi et al..\cite{BFBT}
  
  \section{Linear Solver for the Stokes Equations}
  
We will consider the Stokes equations in the form 
\begin{equation}
    \begin{array}{rrll}
    -\nabla\cdot\left(2\mu\,\varepsilon(\bvec{u})\right) + \nabla p &=& \bvec{f} & \text {in } \Omega \\
    \nabla\cdot \bvec{u} &=& 0 & \text {in } \Omega\\
    \bvec{u} &=& 0 &\text {in } \partial\Omega,
    \end{array}
    \label{eq:stokes}
    \end{equation}
where $\bm{u}$ denotes the fluid's velocity, $p$ it's pressure, and $\mu(\bm{x})$ it's viscosity. Here we use the strain-rate tensor $\varepsilon\left(\bm{u}\right) = \frac12 \left(\nabla \bm{u} + (\nabla\bvec{u})^T\right)$ as we will be considering non-constant viscosity. The right hand side $\bm{f}$ is a forcing term that, inside \aspect{}, typically comes from temperature variation, but for our purposes here, we will only be considering Stokes with velocity and pressure variables and a manufactured right hand side. Finally we are considering incompressible flow with homogeneous Dirichlet boundary conditions, but the developed algorithms also apply to 
the more general case of $\nabla\cdot \bvec{u} = g$, as well as for no-normal-flux boundary constraints
    $
    \nabla \bm{u} \cdot n = 0
    $.

\subsection{Discretization}
Following traditional finite element methods, and using the notation in the work of Kronbichler et al.,\cite{kronbichler_heister_bangerth_2012} we seek coefficients $\bvec{u}_j$ and $p_j$ where, for finite element shape functions $\bm{\varphi^u}_j$ and $\varphi^p_j$,
    \begin{equation}
    \begin{array}{ccc}
    \bvec{u}_h &=& \sum_{j=1}^{N_u} \bvec{u}_j\bm{\varphi^u}_j \\
    p_h &=& \sum_{j=1}^{N_p} p_j\varphi^p_j
    \end{array}
    \end{equation}
    such that the weak formulation
    \begin{equation}
	\begin{array}{rrcl}
		a(\bm{\varphi^u}_i,\,\bm{u}_h) + b(\bm{\varphi^u}_i,\,p_h) & = & f(\bm{\varphi^u}_i) \\
		b(\bvec{u}_h,\,\varphi^p_l)                                  & = & 0,                   
	\end{array}
	\label{eq:weak_form}
\end{equation}
defined by 
$$
    a(\bm{\varphi^u}_i,\,\bm{u}_h) = \int_{\Omega} 2\mu\, \varepsilon(\bm{\varphi^u}_i) : \varepsilon(\bm{u}_h) \;\text{d}x
    \qquad
    b(\bm{\varphi^u}_i,\,p_h) = -\int_{\Omega} (\nabla\cdot \bm{\varphi^u}_i)p_h  \;\text{d}x
    \qquad
    f(\bm{\varphi^u}_i) = \int_\Omega \bm{\varphi^u}_i\cdot\bvec{f}  \;\text{d}x,
    $$
holds for each $1\leq i \leq N_u$ and $1\leq l \leq N_p$ where $N_u$ and $N_p$ are the total number of degrees of freedom for velocity and pressure respectively. We will choose the shape functions from the Taylor-Hood elements $[\QQ_{k+1}]^{\text{dim}}\times\QQ_k,$ $k\geq 1$, which are known to be stable for the system considered here.\cite{guermond_2010,brennerscott02}
    
    Solving this system for coefficients $U = \left\{u_i\right\}$ and $P=\left\{p_j\right\}$ is then equivalent to solving the block linear system
    \begin{equation} 
    \left(\begin{array}{cc} A&B^T\\B&0\end{array}\right) 
    \left(\begin{array}{c}U\\P\end{array}\right)
    =
    \left(\begin{array}{c}F\\0\end{array}\right)
    \label{eq:stokes_linear_system}
    \end{equation}
    where
    $$
    A_{ij} = \int_{\Omega} 2\mu\,\varepsilon(\bm{\varphi^u}_i) : \varepsilon(\bm{\varphi^u}_j) \;\text{d}x
    \qquad
    B_{ij} = -\int_{\Omega} (\nabla\cdot \bm{\varphi^u}_j)\varphi^p_i \;\text{d}x
    \qquad
    F_j = \int_\Omega \bm{\varphi^u}_j\cdot\bvec{f} \;\text{d}x.
    $$

    \subsection{Linear Solver}\label{sec:linear-solver}
    
    There are three common approaches used in the literature for solving \eqref{eq:stokes_linear_system} on large scales: (i) a pressure corrected, Schur complement CG scheme, using multigrid as an approximation to the velocity block,\cite{Gmeiner2015}, (ii) a block-preconditioned Krylov method, also using multigrid on the velocity block, and (iii) an all-at-once multigrid performed on the entire Stokes system, using Uzawa-type smoothers.\cite{Gmeiner2015,BAUER201960}
    For method (ii), there are two main types:
    \begin{enumerate}
     \item[(a)] GMRES\cite{MAY2015496,RudiMalossi2015} (or any Krylov method not requiring symmetry) with block-triangular preconditioner
         \begin{equation}
     P = \left(\begin{array}{cc} A&B^T\\0&-S\end{array}\right),
     \label{eq:stokes_precond}
    \end{equation} 
     or
     \item[(b)] MINRES\cite{Gmeiner2015} with block-diagonal preconditioner
         \begin{equation}
     P_D = \left(\begin{array}{cc} A&0\\0&-S\end{array}\right)
     \label{eq:stokes_precond_diag}
    \end{equation} 
    \end{enumerate}
    where $S = BA^{-1}B^T$ is the \textit{Schur complement}.\cite{elman_wathen_silvester}
    
In the work of Gmeiner et al.,\cite{Gmeiner2015} a comparison is given for each method (using MINRES for (ii)), and the all-at-once multigrid method (iii) was shown to give the fastest and most consistent convergence, as well as the lowest memory overhead. However, there is no comparison that we are aware of for method (ii) using the preconditioner $P$ defined in \eqref{eq:stokes_precond}, and both method have been used separately to demonstrate impressive scaling results (see the work of Bauer et al.\cite{BAUER201960} for method (iii) and Rudi et al.\cite{RudiMalossi2015} for method (ii) with GMRES). 
    
    For the computations in this paper, we will choose method (ii), as this is the framework used in \aspect{}. For preconditioning, we choose $P$ defined in \eqref{eq:stokes_precond}, since, as illustrated by Table~\ref{tab:diag-prec-test}, this leads to roughly half the solve time for a GMRES solve as compared to using $P_D$. This is due to the fact that the block-triangular preconditioner is roughly the same cost to apply as the block-diagonal preconditioner (only containing an extra $B^T$ matrix-vector product), but leads to roughly half the total GMRES iterations.\footnote{It should be noted that the majority of the time $P_D$ is used with MINRES instead of GMRES. However, we do not expect the MINRES to significantly outperform GMRES in runtimes.}   
    
    \begin{table}[tp]
		\centering
		\begin{tabular}{|l|ccc|ccc|}		
		\hline
		Preconditioner:&\multicolumn{3}{c|}{\textbf{$P$}}&\multicolumn{3}{c|}{ \textbf{$P_D$}}\\ 
			\hline
			& iters & Solve time & Time/iter & iters & Solve time & Time/iter\\
			\hline
			860K DoFs & 26 & 0.671 & 0.026 & 54 & 1.294 & 0.020\\
			1.8M DoFs & 27 & 1.829 & 0.068 & 56 & 4.799 & 0.086\\
			3.8M DoFs & 28 & 5.320 & 0.190 & 55 & 10.457 & 0.190\\
			7.7M DoFs & 28 & 11.416 & 0.408 & 57 & 21.142 & 0.371\\
			\hline
		\end{tabular}
		\caption{Iterations required to reduce the residual by $10^6$, time to solution (in seconds), and average time per iteration of a GMRES solve for the Sinker benchmark defined in Section~\ref{sec:nsinker} (4 sinkers and a viscosity ratio of $10^4$), on a series of adaptively refined meshes, and using 32 cores on a workstation . $\hat{A}$ and $\hat{S}$ approximations are defined in Section~\ref{sec:linear-solver}. Here we test the block-triangular preconditioner $P$ against the block-diagonal preconditioner$P_D$. Using $P$ results in roughly half the iterations and solve time compared to $P_D$.}
		\label{tab:diag-prec-test}
\end{table}

    Now, for the preconditioner $P$, consider the following preconditioned system
    \begin{equation}
    \left(\begin{array}{cc} A&B^T\\B&0\end{array}\right)
    \left(\begin{array}{cc} A&B^T\\0&-S\end{array}\right)^{-1} 
    = \left(\begin{array}{cc} A&B^T\\B&0\end{array}\right) 
    \left(\begin{array}{cc} A^{-1}&A^{-1}B^TS^{-1}\\0&-S^{-1}\end{array}\right)
    = \left(\begin{array}{cc} I&0\\BA^{-1}&I\end{array}\right).
    \end{equation}
    This system has only 1 distinct eigenvalue $\lambda=1$ and will converge in at most 2 GMRES iterations assuming an exact preconditioner. Computing exact representations for $A^{-1}$ and $S^{-1}$ is highly impractical, therefore we seek approximate $\hat{A}^{-1}$ and $\hat{S}^{-1}$, respectively.

    \subsubsection{Choosing $\hat{A}^{-1}$}
    \label{sec:Ahat_stokes}
    Since $A$ comes from a vector Laplacian equation, multigrid would appear to be a logical choice given that these methods are widely known to have convergence independent of mesh size $h$ for elliptic boundary value problems.\cite{braesshackbusch83,trottenberg2001} 
    
\begin{figure}[tp]
	\centering
     \includegraphics[width=0.3\textwidth]{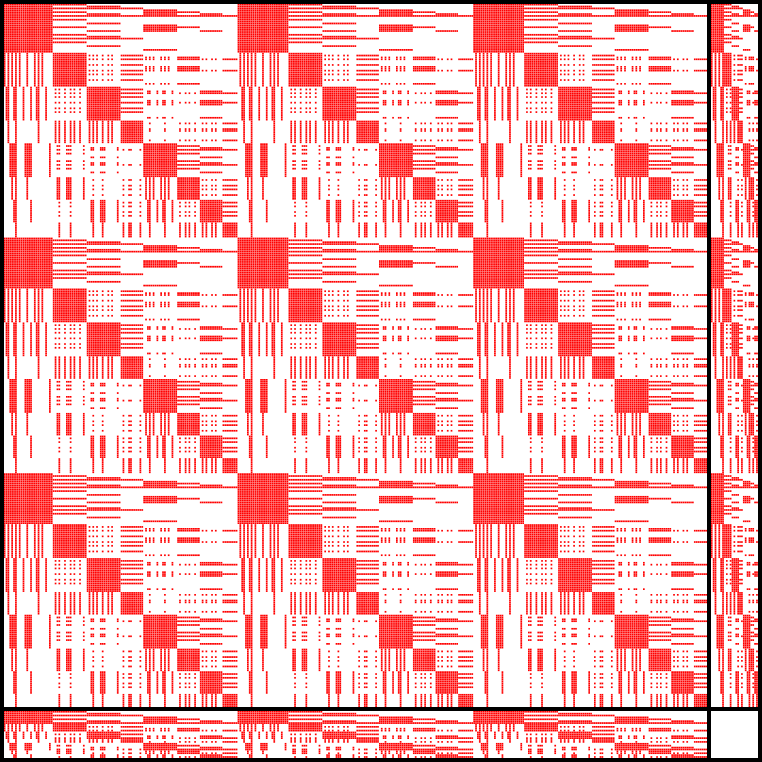}
     \hspace{2.5cm}
     \includegraphics[width=0.3\textwidth]{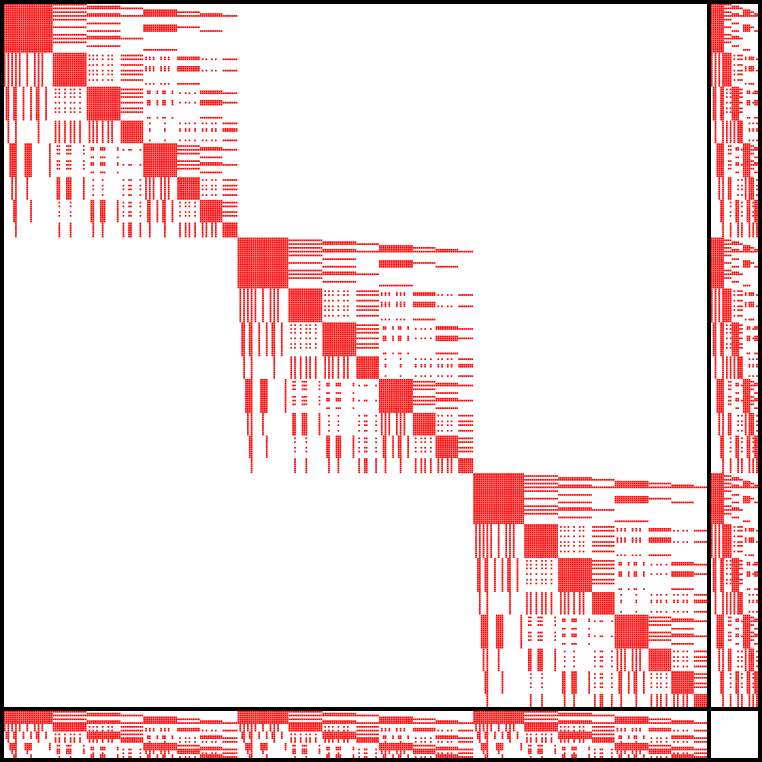}
     \caption{Sparsity pattern for the matrix \eqref{eq:stokes_linear_system} for full coupling of vector components (left) and partial coupling (right). 3D Stokes with 402 degrees of freedom (375 velocity, 27 pressure), ordered component-wise. A red pixel represents a non-zero matrix entry. We see that for partial coupling we have approximately 3x fewer entries as with full coupling.}
     \label{fig:sparsity-stokes}
    \end{figure}

    Currently in \aspect{}, $\hat{A}^{-1}$ is approximated by 1 AMG V-cycle for each Krylov subspace iteration, however the AMG method is not based on the matrix $A$ in \eqref{eq:stokes_linear_system}, but instead we consider the matrix $\hat{A}$ defined as
    \begin{equation}
    \hat{A}_{ij} = \sum_{d=1}^{\text{dim}}\int_{\Omega} 2\mu\, \left[\varepsilon(\bm{\varphi_i^u})\right]_{d,d} : \big[\varepsilon(\bm{\varphi_j^u})\big]_{d,d}
    \label{eq:partial_coupling}
    \end{equation}
    where $\left[\varepsilon(\bm{\varphi_i^u})\right]_{d,d} = \left[\nabla\bm{\varphi_i^u}\right]_{d,d}$ is the $d$th diagonal entry of the gradient operator. The reasons for using this partial coupling of velocity components are the following:
    \begin{enumerate}
	\item[(i)] AMG methods, which depend on the sparsity structure of the underlying matrix, tend to deteriorate when coupling 
	      vector components in higher order computations,\cite{geenen2009} and
	\item[(ii)] the resulting matrix $\hat{A}$ will have far fewer entries (1/3 the entries in 3D, see, e.g., Figure~\ref{fig:sparsity-stokes}), and therefore less storage requirements, faster AMG setup and faster application. 
\end{enumerate}

    We will test this implementation of $\hat{A}^{-1}$ using an AMG v-cycle against one using the GMG v-cycle developed in the work of Clevenger et al.\cite{clevenger_preprint} For this geometric method, we consider a matrix-free v-cycle based on a degree 4 Chebyshev smoother with one smoothing step per level. In contrast with the AMG v-cycle, the level matrices of the GMG v-cycle are based on the fully coupled system $A$, not the partially coupled system $\hat{A}$ since, for the geometric method, there should be no deterioration when using a fully coupled system (contrast to item (i) above), and, as we are using matrix-free operators where each matrix entry is computed on-the-fly at the quadrature level, using the strain rate tensor only consists of adding the off diagonal term in the correct place and dividing by 2; essentially, it is for free (contrast to item (ii) above). From item (iii) above, this could give the GMG-based method an advantage over an AMG-based method in terms of the effectiveness of the preconditioner, while not requiring any extra storage or computational time.

    \subsubsection{Choosing $\hat{S}^{-1}$}\label{sec:S-hat}

    A common choice for approximating $S = BA^{-1}B^T$ is a weighted pressure mass matrix $M_p$, where $\left[M_p\right]_{i,j} = \int_{\Omega}\mu^{-1}\varphi^p_i,\,\varphi^p_j \;\text{d}x$. \cite{Silvester1994,elman_wathen_silvester,kronbichler_heister_bangerth_2012} The reasons for this is that $S$ and $M_p$ are spectrally equivalent for constant viscosity,\cite{elman_wathen_silvester} making $M_p^{-1}$ a good approximation to $S^{-1}$, while $M_p^{-1}$ is far easier to compute. The application of $\hat{S}^{-1}$ then is a simple CG solve with an ILU preconditioner (for matrix-based AMG method) or a Chebyshev iteration (for matrix-free GMG method), converging in between 1-5 iterations. Compared with $\hat{A}^{-1}$, it is not computationally significant.\cite{kronbichler_heister_bangerth_2012} Since the application of $P$ now requires a CG solve whose iteration count may change between applications, a flexible Krylov subspace method must be used for the outer iteration. We will use the flexible variant of GMRES referred to as FGMRES. This solver has the same guaranteed convergence property as GMRES\cite{Saad1993}, however it will require the storage of 2 temporary vectors per iteration (up to the restart length).
    
    It should be noted that this Schur complement approximation begins to break down for large $\text{DR}(\mu)$ (see the work of Rudi et al.\cite{BFBT}) and a more sophisticated approach may be required in those cases. 
    
    \subsection{Viscosity Averaging}
    The viscosity $\mu(\bm{x})$ is evaluated at each quadrature point on the finest cells and then
    averaged using harmonic averaging to obtain a piece-wise constant representation per cell. Since in many applications the viscosity will not come from a functional representation (and therefore cannot be easily evaluated on coarser grids), we must have a way to transfer this cell-wise viscosity on the active mesh to a viscosity on the meshes throughout the level hierarchy. We accomplish this by transferring these coefficients using the same grid restriction operator as discussed in the work of Clevenger et al.,\cite{clevenger_preprint} the work of Brenner and Scott,\cite{brennerscott02} and the work of Janssen and Kanschat,\cite{janssenkanschat11} based on a $\QQ_0^{\text{disc}}$ element (one constant per cell). Essentially, this involves averaging the active cells viscosity, and then recursively setting the the viscosity of a parent cell in the hierarchy to the average of the viscosity of each of its children.
    
    It should be noted that, for problems where the viscosity is represented by a smooth function, this harmonic averaging of the viscosity on the active mesh can result in the loss of convergence order of the solution. However from the work of Heister et al.,\cite{heister_dannberg_gassmoeller_bangerth_2017}, for problems with highly discontinuous viscosity,  averaging viscosity can be beneficial for controlling under- and over-shoots in the numerical approximations of the solution gradient. We do not consider alternatives
    to harmonic averaging here.

\section{Results}\label{sec:results}

  \subsection{Software}
  For all computations we will be using the open-source library \dealii{},\cite{dealII91} which offers 
  scalable parallel algorithms for finite element computations. The \dealii{} library uses functionality from other libraries such as Trilinos\cite{trilinos-overview} (for linear algebra, including Trilinos ML AMG preconditioner) and p4est\cite{bursteddewilcoxghattas11} (for mesh partitioning). 
  
  As mentioned above, we will be comparing our solver to the one used by \aspect{}. \aspect{} is an open-source library written on top of \dealii{} for the specific purpose of solving problems related to the earth's mantle. The method
  presented here is avialable in the current development version of \aspect{}.
  
%
 
\subsection{Benchmark Problem}\label{sec:nsinker}
We give results demonstrating the scalability of the GMG-based method developed here, as well as a comparison with the AMG-based method used in \aspect{}. The test problem used is the ``Sinker'' benchmark described in the work of May et al.\cite{MAY2015496} and the work of Rudi et al..\cite{BFBT} It consists of solving the Stokes problem \eqref{eq:stokes} on the unit cube domain, with $n$ randomly positioned ``sinkers" of higher viscosity throughout the domain. By specifying $\text{DR}(\mu)$, we define a smooth viscosity by $\mu(\bm{x})\in [\mu_{\min},\mu_{\max}]$ where for $X(\bm{x})\in [0,1]$
	\begin{align*}
	X(\bm{x}) &= \prod_{i=1}^{n} \left[1-\exp\left(-\delta\max\left[0,|\bm{c}_i-\bm{x}|-\frac{\omega}{2}\right]^2\right)\right], \\
	\mu(\bm{x}) &= X(\bm{x})\mu_{\min} + (1-X(\bm{x}))\mu_{\max}.
	\end{align*}
	Here $\mu_{\min} = \text{DR}(\mu)^{-1/2}$, $\mu_{\max} = \text{DR}(\mu)^{1/2}$, $\bm{c}_i$ are the center of each sinker, $\delta=200$ controls the exponential decay of the viscosity, and $\omega=0.1$ is the diameter of the sinkers. The right hand side is given by $\bvec{f}(\bm{x}) = (0,0,\beta(X(\bm{x})-1))$ with $\beta=10$ and we use homogeneous Dirichlet boundary conditions for the velocity. Physically, this represents gravity pulling down the high viscosity sinkers. Figure~\ref{fig:sinker_solution} gives a representation of both the velocity and the pressure solution of this benchmark. 
	
   The problem difficulty can be increased by increasing $n$ or $\text{DR}(\mu)$. Table~\ref{tab:sinker_iters_difficulty} gives the iterations required to reduce the residual of the outer FGMRES solve by $10^6$ for different values of these parameters. We see that both AMG and GMG deteriorate as both $n$ and $\text{DR}(\mu)$ increase, with GMG being slightly more robust. This problem can likely be addressed using methods derived in the work of Rudi et al.\cite{BFBT}, where it was shown that this deterioration is due to the approximation loss in the Schur complement solve, and a more sophisticated Schur complement approximation was proposed based on least squares communicators. For the remaining results, we will only consider $n=4$ and $\text{DR}(\mu) = 1e4$, as this is within the range of problems where our Schur complement approach is sufficient.
   
   \begin{table}[tp]
    \begin{minipage}{0.6\linewidth}
    \centering
 \includegraphics[width=1\textwidth]{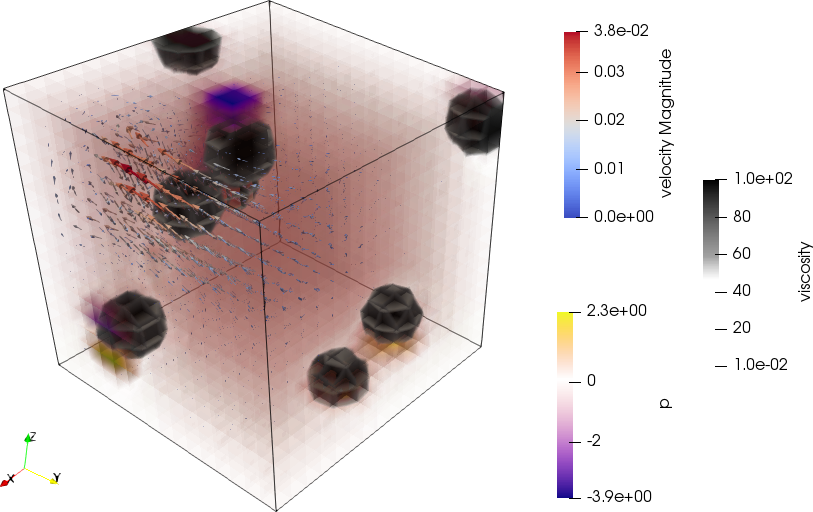}
    \captionof{figure}{$n=8$, $\text{DR}(\mu)=1e4$}
 \label{fig:sinker_solution}
  \end{minipage}
  \hfill
    \centering
    \begin{tabular}{|l|ccc|} 
			\hline
			DR($\mu$) & 1e2 & 1e4 & 1e6 \\  
			\hline\hline
			\textbf{AMG}&&&\\ 
			\hline
			4 sinkers&46&48&52 \\
			8 sinkers&81&196&229 \\
			12 sinkers&76&171&237 \\
			\hline\hline
			\textbf{GMG}&&&\\ 
			\hline
				4 sinkers&20&24&26 \\
				8 sinkers&37&72&128 \\
				12 sinkers&39&80&141 \\
			\hline
		\end{tabular}
  \caption{FGMRES iterations required to reduce the residual by $10^6$ for the sinker benchmark with increasing $n$ and $\text{DR}(\mu)$. Run on a 3D mesh with 860K degrees of freedom ($[\QQ_{2}]^{\text{dim}}\times\QQ_1$ element), distributed over 32 cores.}
		\label{tab:sinker_iters_difficulty}
\end{table}

All timings in this section were from computations run on either Frontera or Stampede2, both at The University of Texas at Austin's Texas Advanced Computing Center. For the Frontera runs, we will be using the Intel Xeon Platinum 8280 (Cascade Lake) nodes which have 56 cores and 192GB per node, and for the Stampede2 runs, we will be using the Intel Xeon Platinum 8160 (Skylake) nodes which have 48 cores and 192GB per node. Both support AVX-512 instructions allowing for vectorization over 8 doubles. The \dealii{} version used is 9.1.0-pre, and we compile using gcc 7.1.0, intel-mpi 17.0.3. The p4est version is 2.0.0, the Trilinos version is 12.10.1 and the \aspect{} version is 2.1.0-pre. 
   
Each line in the timing plots connect the median time of 5 distinct runs, with all 5 runs shown as points. We will test strong scaling (problem size stays constant, number of cores increase) and weak scaling (problem size per core stays constant) for key parts of the computation. Tests run on adaptively refined meshes will include the model for partition imbalance described in the work of Clevenger et al..\cite{clevenger_preprint} This model should represent the ideal scaling we can expect to see for the GMG method given the imbalance of cells in the level hierarchy.

\subsection{AMG/GMG Comparison}\label{sec:amg_gmg_comp}
For a comparison between the AMG and GMG preconditioners, we will consider an adaptively refined mesh, where for each refinement, the number of cells are roughly doubled. We start with a mesh of 4 global refinements and create each new mesh using a gradient based estimator to refine roughly 1/7 of the cells from the previous refinement cycle, doubling the number of cells in our mesh. 
Trilinos ML AMG is used with a Chebyshev smoother containing two sweeps and an aggregation threshold of 0.001. These values are chosen here as they are the default values in \aspect{} and represent the values used by the majority of the users of the code. 

Table~\ref{tab:amg_gmg_48} gives the runtimes for such a mesh with 5 levels of adaptive refinement on 48 cores ($18.5\times 10^6$ degrees of freedom). The ``Setup'' time includes the distribution of the degrees of freedom, setting up of any sparsity patterns necessary, as well as the setup of the data structures required for the matrix-free GMG transfer. Here, our GMG method requires roughly 2x the work for distributing the degrees of freedom, but does not need to build any sparsity patterns. This results in roughly equivalent setup times between AMG and GMG, with GMG being slightly faster. In theory, we should easily be able to lower the requirement of 2x the work in degree of freedom distribution. The ``Assemble'' timing includes all matrix assembly (system matrix, preconditioner matrix, AMG setup) as well as assembling the right hand side of the linear system and vectors/tables related to the matrix-free operators. Here is where we see the largest advantage for the GMG method as it has no matrices to assemble, resulting in more than a 10x faster assembly. Combining setup and assembly with the linear solve, we have that the GMG method is around 3x faster for this problem. 

\begin{table}[tp]
    \centering
    \begin{tabular}{|l|rr|r|}
 \hline
 &AMG&GMG&factor\\
 \hline
 Setup&12.6s&10.3s&1.2x\\
 Assemble&32.5s&2.9s&11.2x\\
 Solve&38.6s&14.8s&2.6x\\
 \hline
 Total&83.7s&28.0s&3.0x\\
\hline
\end{tabular}
  \caption{Timing comparison between AMG and GMG for an adaptively refined, 3D mesh, with $18.5\times10^6$ DoFs ($[\QQ_{2}]^{\text{dim}}\times\QQ_1$ element) on one node (48 cores). Timings are from the Stampede2 machine.}
 \label{tab:amg_gmg_48}
\end{table}

However, for time dependent applications, many time steps will typically be solved without further refining the mesh, in which case, we no longer need to call the ``Setup'' functionality. Then the program time will be dominated by assembly and solve, in which case GMG will be about 4x faster here.

\begin{figure}[tp]
\centering
		\subfigure[One application of the Stokes block preconditioner.]{
		\includegraphics[width=0.47\textwidth]{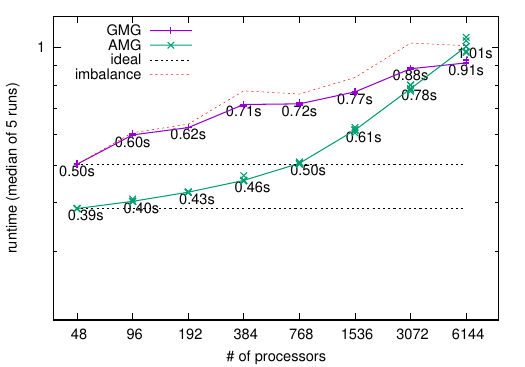}
		\label{fig:amg_gmg_prec}
		}
		
		\subfigure[FGMRES solve.]
		{
		\includegraphics[width=0.47\textwidth]{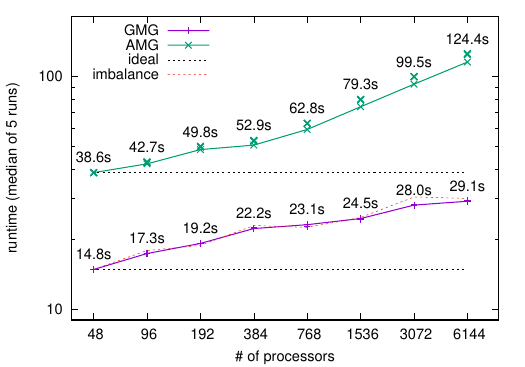}
		\label{fig:amg_gmg_solve}
		}
		\subfigure[FGMRES solver speedup (DoFs/time). Labels represent weak scaling efficiency, parentheses is efficiency to the partition imbalance.]
		{
		\includegraphics[width=0.47\textwidth]{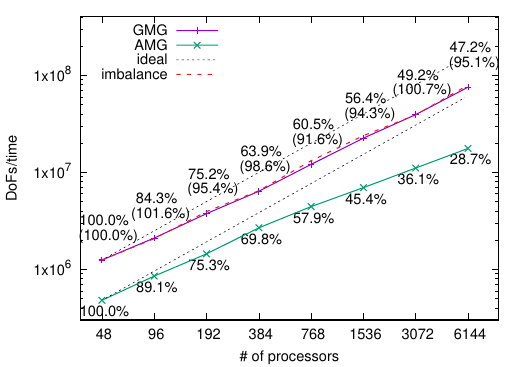}
		\label{fig:amg_gmg_solve_speedup}
		}
		\caption{Weak scaling comparison, adaptive refinement, 3D mesh, $[\QQ_{2}]^{\text{dim}}\times\QQ_1$ element. Timings are from the Stampede2 machine.}
	\label{fig:weak}
\end{figure}

Expanding on this, we look at the weak scaling of each component up to 6,144 cores and mesh size of $2.2\times10^9$ degrees of freedom. Starting with the linear solve, as with the scaling plots for global refinement above, Figure~\ref{fig:amg_gmg_prec} gives the timing for the preconditioner application, and Table~\ref{tab:amg_gmg_iters} gives the number of FGMRES iteration required in the solve. Here we see that, while the AMG preconditioner is cheaper to apply for all but the last data point, the iteration counts for GMG are much lower and stay constant while the AMG iteration counts increase by over $50\%$. Figure~\ref{fig:amg_gmg_solve} shows the solve time and Figure~\ref{fig:amg_gmg_solve_speedup} shows the speedup. The red dashed line in this plot represents the ideal weak scaling (black dashed line) multiplied by a factor representing the imbalance of the parallel mesh partition that occurs on computations using adaptive mesh refinement, explain in the work of Clevenger et al..\cite{clevenger_preprint} This factor represents the increase in runtime one can expect with the current mesh partition (as opposed to a fair distribution), and is both dependent on the mesh refinement scheme and the number of cores used.  Here we see that, even with the imbalance of the partition, the scaling for GMG is more efficient than AMG, and there is near perfect efficiency when taking into account the imbalance of the mesh partition, which, according to the work of Clevenger et al.,\cite{clevenger_preprint} should remain bounded.

 \begin{table}[tp]
    \centering
    \begin{tabular}{|r|r|cc|}
		\hline
			Procs & DoFs& AMG & GMG  \\  
			\hline
			48&18M&53&27 \\
			96&36M&56&27 \\
			192&72M&62&28 \\
                        384&141M&62&28 \\
			768&278M&68&28 \\
			1536&551M&75&28 \\
       			3072&1.1B&80&28 \\
			6144&2.2B&83&28 \\
			\hline
		\end{tabular}
  \caption{FGMRES iterations required to reduce the residual by $10^6$ on problems depicted in Figure~\ref{fig:amg_gmg_prec}.}
		\label{tab:amg_gmg_iters}
\end{table}

\begin{figure}[tp]
\centering
		\subfigure[Setup.]{
		\includegraphics[width=0.47\textwidth]{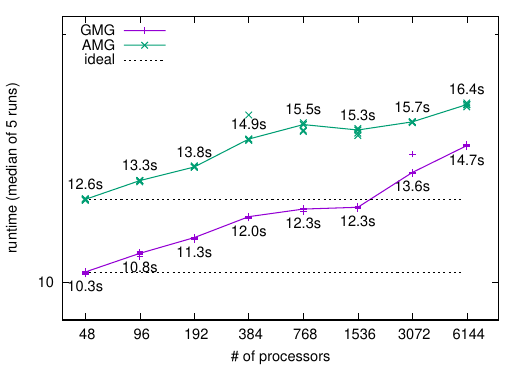}
		\label{fig:amg_gmg_setup}
		}
		\subfigure[Assembly.]
		{
		\includegraphics[width=0.47\textwidth]{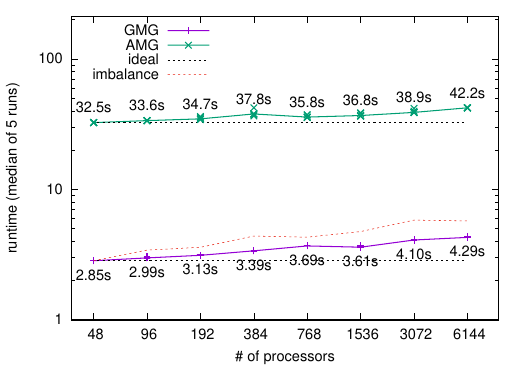}
		\label{fig:amg_gmg_assemble}
		}
		\caption{Setup/Assembly of the linear system. Timings are from the Stampede2 machine.} 
\end{figure}

The weak scaling of the setup is shown in Figure~\ref{fig:amg_gmg_setup}, and again we see that AMG and GMG are roughly equivalent, with GMG being slightly faster. The setup of the linear system should be optimized in the future since we are roughly on the same order of magnitude as the solve time, and efforts should be made to improve the scaling. The weak scaling of the assembly is shown in Figure~\ref{fig:amg_gmg_assemble}. Unsurprisingly, the GMG assembly (consisting mostly of viscosity evaluation) is much cheaper than AMG as there are no matrices to assemble.

\subsection{Memory Consumption}\label{sec:memory-consumption}

The storage requirements for FGMRES, required by the preconditioner in \aspect{} due to the non-constant nature of the CG solve for the mass matrix (discussed in Section~\ref{sec:S-hat}), is quite restrictive as one will need to store two additional vectors for each Krylov subspace iteration up to the restart length. With a restart length of 50, as is the default inside \aspect{}, one will have to allocate the storage of up to 101 vectors inside the Krylov solve alone. This can greatly affect the overall storage of the GMG-based method, as seen in Table~\ref{tab:amg_gmg_memory}.

Table~\ref{tab:amg_gmg_memory} gives an estimation of the memory consumptions (in MB) of the largest global objects required by each method on a globally refined mesh with 113K degrees of freedom run on a single core. While, for the AMG-based method the temporary vectors make up around 35\% of the total storage, these vectors represent around 95\% of the storage required by the GMG-based method. Therefore, we only see around 2.7x less memory used for the GMG-based method, even though there are no matrices that need to be stored. 

\begin{table}[tp]
	\centering
	\begin{tabular}{|l|rr|}
		\hline
		Memory (MB)   & AMG   & GMG  \\
		\hline
		Triangulation & 1.9   & 1.9  \\
		DoFHandlers   & 2.8   & 5.7  \\
		Constraints   & 1.0   & 2.7  \\
		$A$           & 174.2 & -    \\
		$B$ and $B^T$ & 31.4  & -    \\
		$\hat{A}$     & 58.5  & -    \\
		$\hat{S}$     & 1.4   & -    \\
		Vectors(105)   & 179.6  & 179.6 \\
		AMG matrices  & 59.8  & -    \\
		\hline
		Total         & 510.6 & 189.9 \\
		\hline
	\end{tabular}
	\caption{Memory consumption required for major components of AMG- and GMG-based methods for globally refined, 3D mesh, with 113K DoFs ($[\QQ_{2}]^{\text{dim}}\times\QQ_1$ element) on 1 core. The total number of vectors required is set at 105. This estimates takes into account the temporary vectors in the FGMRES solve as well as others required throughout the code (right-hand side vector, solution vector, and two previous timesteps are common).}
	\label{tab:amg_gmg_memory}
\end{table}

The first step to remove this requirement for temporary vectors is to find a preconditioner which does not require a flexible Krylov method. This can easily be done by considering a further approximation of the Schur complement $S$ by taking $\hat{S} = \text{diag}(M_p)$, as is done in the work of May et al.\cite{MAY2015496} Another option is to approximate $S$ as one multigrid v-cycle of the mass matrix $M_p$, which we will do for the remainder of the results. As with the CG solve of the mass matrix (described in Section~\ref{sec:S-hat}), the application of the vcycle will not be computationally significant as compared to the application of $\hat{A}^{-1}$. With this new Schur complement approximation, one can use GMRES instead of FGMRES, which will reduce the number of temporary vectors by half. Then, for the test in Table~\ref{tab:amg_gmg_memory}, we now have a total of 94.1MB of storage for vectors out of a total of 104.4MB for the GMG-based method (90\% of the total), and we see around 4.9x less memory as opposed to the current AMG-based method in \aspect{}. 

While this is a step in the right direction, we would like a method whose storage was constant with respect to the number of iterations. We will seek further improvements using the IDR($s$) method from the work of Gijzen and Sonneveld.\cite{Gijzen2008} These Krylov methods have a short-term recurrence, requiring the storage of $5+3s$ vectors for the solve, and there are $s+1$ matrix-vector products and preconditioner applications per iteration. Assuming the correct choice of parameters, IDR(1) is equivalent to BiCGStab, and, as $s$ is increased, the convergence rate has been demonstrated to approach the convergence rate of GMRES.\cite{Gijzen2008}  The implementation in \dealii{}is based on the preconditioned IDR($s$) from the work of Gijzen and Sonneveld\cite{elegant-idr}. We consider IDR(2) for the computations in this section. 

Table~\ref{tab:idr-gmres} gives a comparison between the strong scaling and storage requirements of GMRES and IDR(2) for the Sinker benchmark run on a globally refined mesh with $27\times10^9$ unknowns. Since the preconditioner application dominates the overall runtime of the solve, it is no surprise that the IDR(2) method is slightly more expensive. However, the times are comparable and both methods see similar scaling. The IDR(2) method, on the other hand, requires roughly 1/3 of the total vector storage.

\begin{table}[tp]
		\centering
		\begin{tabular}{|r|cccc|cccc|}		
		\hline
		&\multicolumn{4}{c|}{GMRES}&\multicolumn{4}{c|}{IDR(2)}\\ 
			\hline
			Cores & iters & solve time& \# of $P^{-1}$ & vectors & iters & solve time & \# of $P^{-1}$ & Vectors\\
			\hline
			14,336 & 30 & 40.5 & 30 & 31 & 11 & 42.1 & 33 & 11\\
			28,672 & 30 & 23.1 & 30 & 31 & 11 & 21.8 & 33 & 11\\
			57,344 & 30 & 13.0 & 30 & 31 & 13 & 13.5 & 39 & 11\\
			114,688 & 30 & 7.74 & 30 & 31 & 13 & 7.90 & 39 & 11\\
			\hline
		\end{tabular}
		\caption{Iterations required to reduce the residual by $10^6$, time to solution (in seconds), number of preconditioner applications, and vectors stored for a solve of the Sinker benchmark defined in Section~\ref{sec:nsinker} (4 sinkers and a viscosity ratio of $10^4$), on a globally refined mesh with $27\times10^9$ DoFs. Timings are from the Frontera machine.}
		\label{tab:idr-gmres}
\end{table}

Returning then to the test in Table~\ref{tab:amg_gmg_memory}, consider Table~\ref{tab:all_memory} which gives the storage requirement for all FGMRES(50), GMRES(50), and IDR(2).
We see that, for the GMG-based method, switching from GMRES to FGMRES leads to a 1.8x memory increase, and switching from IDR(2) to FGMRES leads to a 5.3x memory increase. If we compare the increase in memory from IDR(2) with the GMG-based method to FGMRES with the AMG-based method, we can expect a roughly 14x memory increase. This implies that the capabilities of solving large scale problems inside \aspect{} would be drastically increased by using the GMG-based Stokes solve with IDR(2), allowing for runs with an order of magnitude more unknowns. 

\begin{table}[tp]
	\centering
	\begin{tabular}{|l|cc|c|c|}
		\hline
		& \multicolumn{2}{c|}{FGMRES(50)} & GMRES(50) & IDR(2) \\ 
		\hline
		Memory (MB)   & AMG   & GMG   & GMG   & GMG  \\
		\hline
		Triangulation & 1.9   & 1.9   & 1.9   & 1.9 \\
		DoFHandlers   & 2.8   & 5.7   & 5.7   & 5.7 \\
		Constraints   & 1.0   & 2.7   & 2.7   & 2.7 \\
		$A$           & 174.2 & -     & - & - \\
		$B$ and $B^T$ & 31.4  & -     & -  & - \\
		$\hat{A}$     & 58.5  & -     & -  & - \\
		$\hat{S}$     & 1.4   & -     & -   & - \\
		Vectors       & 179.6 & 179.6 & 94.1 & 25.7 \\
		AMG matrices  & 59.8  & -     & -  & -  \\
		\hline
		Total         & 510.6 & 189.9 & 104.4 & 36.0 \\
		\hline
	\end{tabular}
	\caption{Memory consumption required for major components of AMG- and GMG-based methods for globally refined, 3D mesh, with 113K DoFs ($[\QQ_{2}]^{\text{dim}}\times\QQ_1$ element) on 1 core. FGMRES(50) (105 total vectors) estimates are given for both AMG and GMG, and both GMRES(50) (55 total vectors) and IDR(2) (15 total vectors) estimates are given for GMG.}
	\label{tab:all_memory}
\end{table}

\subsection{Large Scale Computations}\label{sec:idr-scaling}

\begin{figure}[tp]
\centering
		\subfigure[One application of the Stokes block preconditioner.]{
		\includegraphics[width=0.47\textwidth]{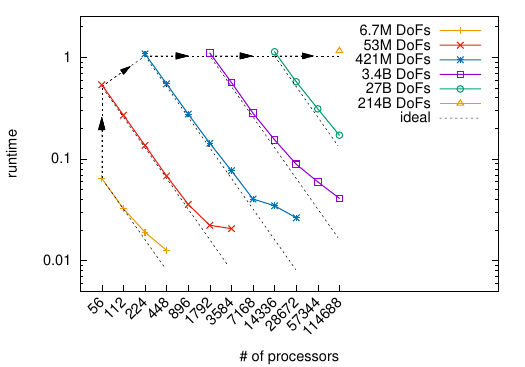}
		\label{fig:gmg_prec_strong}
		}
		\subfigure[IDR(2) solve.]
		{
		\includegraphics[width=0.47\textwidth]{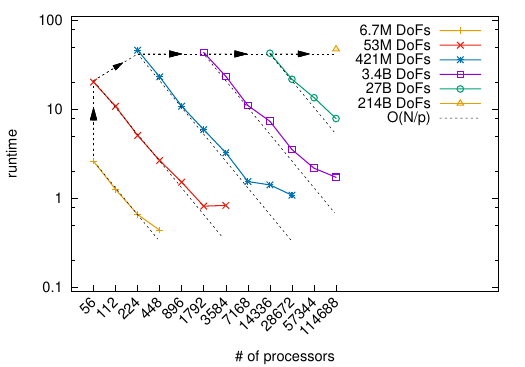}
		\label{fig:gmg_solve_strong}
		}
		\caption{Strong and weak scaling for GMG-based method for globally refined, 3D mesh with $[\QQ_{2}]^{\text{dim}}\times\QQ_1$ element. Timings are from the Frontera machine.} 
	\label{fig:gmg_strong}
\end{figure}

Figure~\ref{fig:gmg_prec_strong} gives the strong scaling for an application of the Stokes preconditioner from a globally refined mesh with between 6-11 refinement levels, and Table~\ref{tab:gmg_iters_strong} gives the corresponding IDR(2) iteration counts to reduce the residual by $10^6$. The dashed scaling lines here are all computed based on the data point at 56 cores/$6.7\times10^6$ DoFs, and therefore the plot represents both strong and weak scaling. We see scaling to around 30-60K DoFs/core for the preconditioner and roughly constant iteration counts. Figure~\ref{fig:gmg_solve_strong} gives the timings for the IDR(2) solve, and unsurprisingly, since the iteration counts are almost constant, we see the same scaling. 


 \begin{table}[tp]
    \centering
       \begin{tabular}{|r|rrrrrr|}
	\hline
	Procs$\backslash$DoFs & 6.7M & 53M & 412M & 3.4B & 27B & 217B \\ 
	\hline
	56& 12 & 11 &&&&\\
	112& 12 & 12 &&&&\\
	224& 11 & 11 & 13 &&&\\
	448& 11 & 12 & 13 &&&\\
	896&& 13 & 12 &&&\\
	1792&& 11 & 12 & 12 &&\\
	3584&& 11 & 12 & 12 &&\\
	7168&&& 12 & 12 &&\\
	14336&&& 12 & 14 & 11 &\\
	28672&&& 13 & 12 & 11 &\\
	57344&&&& 12 & 13 &\\
	114688&&&& 12 & 13 & 12\\  	
	\hline
   \end{tabular}
  \caption{IDR(2) iterations required to reduce the residual by $10^6$ for the GMG-based preconditioner on problems depicted in Figure~\ref{fig:gmg_prec_strong} (Globally refined mesh).}
		\label{tab:gmg_iters_strong}
\end{table}

Lastly, we will compare the performance of the GMG v-cycle on the velocity block (representing approximately 70\% of the preconditioner runtime) to the peak performance possible on the Frontera machine. Using the tool LIKWID,\cite{likwid} we find that we must perform approximately 1,575 Flops/DoF for a globally refined mesh. If we consider that the v-cycle application on the velocity space for the problem containing $217\times10^9$ Stokes DoFs ($206\times10^9$ velocity DoFs) takes 0.827s of compute time, we have an arithmetic throughput of 0.40 Petaflops/s. The Intel Xeon Platinum 8280 (Cascade Lake) nodes on Frontera can execute 16 double precision floating point operations per cycle, with a clock speed of 2.7 GHz, and so on 114,688 cores we have a peak performance of 5.0 Petaflops/s. This means that we are operating at 8\% of peak performance. Similar performance for GMG v-cycles using \dealii{} was seen in the work of Arndt et al.,\cite{arndt2019dealii} where a GMG preconditioner based on a scalar Laplacian discretized with a $\QQ_4^{\text{disc}}$ elements was calculated to be around 13\% peak performance. There, they found the computations were memory bound,\cite{arndt2019dealii} which is also expected for the computations here.

\section{Conclusion}
In this article we presented a geometric multigrid method for preconditioning of the velocity block in a Stokes solve with variable viscosity. The presented method was designed using a matrix-free framework and has capabilities to be run on adaptively refined meshes and in parallel. The parallel scalability of this method was shown for both strong and weak scaling, for both global and adaptive refinement, with up to 114,688 cores and up to 217 billion degrees of freedom. We performed a comparison of the developed method with the AMG-based preconditioner currently used in the \aspect{} code and found that the GMG-based preconditioner was both more robust (lower iteration counts and constant convergence with mesh refinement) and exhibits better weak scaling than the AMG method, resulting in at least 3x faster computations than the AMG-based method. The GMG-based method was also estimated to require 14x less memory. All-in-all, the GMG-based Stokes solver in this paper (and implemented in \aspect{}) has allowed for the solution of problems 2 orders of magnitude larger than the largest problem ever solved by the AMG-based method in \aspect{} (these largest runs are given by the 2.2 billion degrees of freedom run in Section~\ref{sec:amg_gmg_comp}). Finally, the GMG-based method, using IDR(2) as a solver, was shown to scaling to around 30K DoFs/core, and the GMG v-cycle for the velocity block, which represents around 70\% of the total preconditioner application, was shown to have a throughput of 8\% of peak performance on 114,688 cores.

\section*{Acknowledgments}

The authors were partially supported by NSF Award EAR-1925575, OAC-2015848,
DMS-2028346, and by the Computational Infrastructure in Geodynamics initiative
(CIG), through the NSF under Award EAR-0949446 and EAR-1550901 and The
University of California -- Davis.

This work used the Extreme Science and Engineering Discovery Environment (XSEDE), which is supported by National Science Foundation grant number ACI-1548562.\cite{XSEDE} The authors acknowledge the Texas Advanced Computing Center (TACC) at The University of Texas at Austin for providing access to both the Stampede2 and Frontera machines that have contributed to the research results reported within this paper. Clemson University is acknowledged for generous allotment of compute time on Palmetto cluster. Lastly, we would like to thank Martin Kronbichler of the Technical University of Munich (TUM) for his extensive help in the understanding and implementation of matrix-free algorithms. 

\bibliographystyle{plain}
\bibliography{bibliography}

\end{document}